\documentclass[11pt,sumlimits,intlimits]{amsart}
\usepackage{url,cite,hyperref,amsmath,amsthm,amssymb,mathtools}
\usepackage[margin=1.25in]{geometry}
\usepackage{xcolor}
\usepackage{comment}
\usepackage[all]{xy}

\newtheorem{thmx}{Theorem}

\newtheorem{theorem}{Theorem}[section]
\newtheorem{proposition}[theorem]{Proposition}
\newtheorem{corollary}[theorem]{Corollary}
\newtheorem{lemma}[theorem]{Lemma}
\theoremstyle{definition}
\newtheorem{definition}[theorem]{Definition}
\newtheorem{remark}[theorem]{Remark}
\newtheorem{example}[theorem]{Example}

\newcommand{\Z}{\mathbb{Z}}
\newcommand{\Q}{\mathbb{Q}}
\newcommand{\id}{\textup{id}}

\newcommand{\C}{\mathbb{C}}
\newcommand{\N}{\mathbb{N}}

\newcommand{\T}{\mathbb{T}}

\newcommand{\ve}{\varepsilon}

\newcommand{\dimnuc}{\dim_{\mathrm{nuc}}}

\numberwithin{equation}{section}

\title[Strong quasidiagonality and decomposition rank]
{Strongly quasidiagonal C*-algebras with nuclear dimension one
and infinite decomposition rank}
\author{Caleb Eckhardt}
\dedicatory{Dedicated to the memory of my mentor, Zhong-Jin Ruan}
\address{Department of Mathematics, Miami University, Oxford, OH, 45056}
\email{eckharc@miamioh.edu}
\begin{document}
\thanks{This work was partially supported by an AMS-Simons research enhancement grant for PUI faculty GR000765.}

\begin{abstract} We construct the first, to our knowledge,  examples of strongly quasidiagonal C*-algebras with finite nuclear dimension and infinite decomposition rank. The examples are extensions of the form $0\to I \to E\to C(\T)\to 0$ where $I$ is a simple AF-algebra. These examples are in some sense the easiest `natural' choices in light of Spielberg's work on AF-embedding extensions and the behavior of each dimension theory with respect to extensions.  The standard extension bound provides $\dimnuc E \leq 2$ which is improved to 1 by recent work of Evington, Faurot and Schafhauser. Nuclear dimension 1 is the lowest possible for such an example because decomposition rank and nuclear dimension agree at 0. Finally, failure of the comparison property used to show infinite decomposition rank was inspired by local operator space techniques of Junge, Ozawa and Ruan.
\end{abstract}
\maketitle

\section{Introduction}
Decomposition rank and nuclear dimension are closely related dimension theories for nuclear C*-algebras. The starkest difference between them may be that finite decomposition rank implies strong quasidiagonality \cite[Theorem 5.3]{Kirchberg04} while finite nuclear dimension does not \cite[Example 6.3]{Winter10}, in fact classifiable purely infinite C*-algebras  have finite nuclear dimension \cite[Theorem 7.5]{Winter10}. There are large classes of C*-algebras for which decomposition rank and nuclear dimension agree and are finite. In all of these cases the examples are strongly quasidiagonal, and until now we had no examples of strongly quasidiagonal C*-algebras whose decomposition rank and nuclear dimension differed.  Let me briefly describe the examples where the dimension values agree. 

Let $A$ be a C*-algebra and let $\operatorname{dr}A, \dimnuc A$ be the decomposition rank, respectively nuclear dimension of $A.$ We first note the trivial bound of $\dimnuc A \leq \operatorname{dr}A.$ 
By \cite[Example 4.1]{Kirchberg04}  we have $A$ is an AF-algebra if and only if $\dimnuc A= \operatorname{dr}A=0.$ From this one also obtains equality if $\operatorname{dr}A \leq 1.$ For example if $A$ is an $A\T$-algebra or more generally an inductive limit of C*-algebras with decomposition rank bounded by 1.
It follows from Winter's article \cite{Winter04} (see also \cite[Corollary 3.21]{Elliott20}) that if $A$ is subhomogeneous then $\dimnuc A= \operatorname{dr} A = \max_k \dim(\operatorname{Prim}_k(A)).$ By \cite[Corollary 2.10]{Winter10}  if $A$ is continuous trace, then $\dimnuc A = \operatorname{dr} A =\dim(\operatorname{Prim} A).$
Finally, if $A$ is classifiable then quasidiagonality implies equality.  More specifically if $A$ is a unital simple separable quasidiagonal nuclear C*-algebra with finite nuclear dimension that satisfies the UCT then the Tikuisis-White-Winter theorem \cite{Tikuisis17} and  \cite[Theorem B]{Castillejos21} show that $\operatorname{dr}A\leq 1$ and therefore $\operatorname{dr} A =\dimnuc A.$

In this paper we give the first examples of strongly quasidiagonal C*-algebras with different  nuclear dimension and decomposition rank.
\begin{thmx}[see Corollary \ref{cor:main}] There are continuum many non-isomorphic strongly quasidiagonal C*-algebras with nuclear dimension 1 and  infinite decomposition rank.
\end{thmx}

Nuclear dimension and decomposition rank behave differently with respect to extensions $0\to I\to E\to B\to 0.$  If $I$ and $B$ both have finite nuclear dimension then so does $E$ by \cite[Proposition 2.9]{Winter10}.  On the other hand if $I$ and $B$ both have finite decomposition rank then $E$ may not have finite decomposition rank. But if the extension is quasidiagonal, i.e. $I$ contains an approximate quasicentral unit of projections then $E$ does have finite decomposition rank \cite[Proposition 6.1]{Kirchberg04}.

Hence one looks for non-quasidiagonal extensions $0\to I\to E\to B\to 0$ where both $I$ and $B$ have finite nuclear dimension and $E$ is strongly quasidiagonal.  
My starting point was extensions reminiscent of my favorite quasidiagonal  C*-algebra coming from a non-quasidiagonal extension: $0\to \mathcal{K}\oplus \mathcal{K} \to C^*(S\oplus S^*)\to C(\T)\to 0$
where $S$ is the unilateral shift. Spielberg's  \cite[Theorem 1.15]{Spielberg88}  provides a relatively easy to check K-theoretic condition to decide (strong) quasidiagonality of extensions $0\to I\to E\to C(\T)\to 0$, where $I$ is a simple AF-algebra (although his theorem is more general than this).  Combining the idea of the above mentioned Toeplitz like extension with some K-theoretic calculations and Spielberg's theorem led directly to the candidates.

Our strategy to show infinite decomposition rank came from Junge, Ozawa and Ruan's \cite[Theorem 3.4]{Junge03}.  
Let $A$ be a finite-dimensional C*-algebra and $x\in A$,  then $x^*x$ is unitarily equivalent to $xx^*.$ They use this simple fact to deduce interesting information about the local operator space structure of nuclear non-finite C*-algebras. We take the same fact about finite dimensional C*-algebras but run a modified version of it through the definition of decomposition rank to obtain a comparability type condition we call $(d,f)$-comparison (see Definition \ref{def:CPdcomp}). 

Returning to the above extension, let $T=S\oplus S^*.$ Even though $C^*(T)$ is stably finite, notice that the projections $1-T^*T$ and $1-TT^*$ are incomparable (Definition \ref{def:CPdcomp}), i.e. for fixed $d$, there is no collection of elements $x_0, ..., x_d \in C^*(T)$ such that $\sum x_ix_i^*$ is arbitrarily close to $1-T^*T$ while $\sum x_i^*x_i$ is arbitrarily close to $1-TT^*.$  This fact shows $C^*(T)$ fails $(d, 1-t)$-comparison  and hence provides an alternate proof that $C^*(T)$ has infinite decomposition rank.  We show that our extensions also violate this comparability condition but, in contrast to $C^*(T),$ are strongly quasidiagonal.

\medskip
 
\noindent\textbf{AI Statement:} A couple years ago I tried to show the extensions from Definition \ref{def:basicext} and similar extensions violated $(d, 1-t)$-comparison and then use Proposition \ref{prop:basicdr} to obtain counterexamples. I was unable to do it.  A discussion with ChatGPT led to Lemma \ref{lem:infdr} which was the missing piece.  I wrote the article.

\section{Results} We start with a definition reminiscent of Cuntz-Pedersen equivalence \cite{Cuntz79}.
\begin{definition}\label{def:CPdcomp} Let $A$ be a C*-algebra, $d\geq 0$ an integer and $f:[0,1]\to [0,1]$ a continuous function.  We say $A$ has $(d, f)$-\textbf{comparison} if for each contraction $x\in A$  there are $d+1$ sequences of contractions $z_{0, m}, ...,z_{d, m}\in A$ (or $\tilde{A}$ if $A$ is not unital) such that
\begin{equation*}
\lim_{m\to \infty} \sum_{i=0}^d z_{i,m}^*z_{i,m} = f(x^*x), \qquad \text{ and }\qquad \lim_{m\to \infty}\sum_{i=0}^d z_{i,m}z_{i,m}^* = f(xx^*),
\end{equation*}
where the functional calculus takes place in the unitization if $A$ is non-unital.
\end{definition}
Trivially all C*-algebras have $(0, \id)$-comparison.  We record a few,  less trivial observations about $(d, f)$-comparison to relate it to more familiar properties. In particular we show $(d,1-t)$-comparison sits between stable rank 1 and (strong) finiteness.
\begin{proposition}\label{prop:sr1} Let $A$ be a unital C*-algebra such that $x^*x$ and $xx^*$ are approximately unitarily equivalent for every $x\in A.$  Then $A$ has $(0, f)$-comparison for all continuous $f:[0, 1]\to [0, 1].$ In particular all C*-algebras with stable rank 1 have $(0, f)$-comparison for all continuous $f:[0, 1]\to [0, 1].$
\end{proposition}
\begin{proof} Let $f:[0,1]\to [0, 1]$ be continuous and let $x\in A$ be a contraction. By assumption there is a sequence $(u_m)$ of unitaries such that $u_m^*x^*xu_m\to xx^*.$  Set $z_{0, m} = u_m^*f(x^*x)^\frac{1}{2}.$ Then $z_{0, m}^*z_{0, m} = f(x^*x)$ and $z_{0, m}z_{0, m}^*\to f(xx^*).$ 

If $A$ has stable rank 1 and $x\in A$, then obtain a sequence $y_m$ of invertible elements that converge to $x.$ By considering the polar decomposition of $y_m$ we have $y_m^*y_m$ is unitarily equivalent to $y_my_m^*$ hence $x^*x$ is approximately unitarily equivalent to $xx^*.$
\end{proof}

\begin{lemma} \label{lem:quotient} Let $f:[0, 1]\to [0, 1]$ be continuous  and $d\geq0.$ Suppose $A$ is a C*-algebra with $(d, f)$-comparison. Then every quotient of $A$ also has $(d, f)$-comparison. 
\end{lemma}
\begin{proof} Let $I$ be an ideal of $A$ and $x\in A/I$ be a contraction.  Then lift $x$ to a contraction  $\tilde{x}\in A$, obtain the $z$ sequences and then quotient back down to $A/I.$  
\end{proof}
\begin{proposition} \label{prop:stabfin} Let $f(t)=1-t$ and suppose that $A$ is a unital C*-algebra with an ideal $I$ such that $A/I$
is not a finite C*-algebra.  Then $A$ does not have $(d, f)$-comparison for any $d\geq0.$
\end{proposition}
\begin{proof} By Lemma \ref{lem:quotient} we may assume that $A$ is not a finite C*-algebra.  Then there is an $x\in A$ such that $1_A-x^*x=0$ and $1_A-xx^*\neq 0.$ It follows that $A$ can not have $(d, f)$-comparison for any $d.$ 
\end{proof}

\begin{proposition}\label{prop:basicdr} Let $A$ be a unital C*-algebra with $\operatorname{dr}A \leq d <\infty.$ Then $A$ has $(d,f)$-comparison for every continuous $f:[0, 1]\to [0,1].$
\end{proposition}
\begin{proof} Let $\ve>0$ and $f:[0, 1]\to [0, 1]$ be continuous. Choose $\delta >0$  so if $0\leq a, b\leq 1$ satisfy $\| a-b \|<\delta, $ then $\| f(a)-f(b) \|<\ve.$  Let $x\in A$ be a contraction.  By \cite[Proposition 5.1]{Kirchberg04} there are finite dimensional C*-algebras $F_0, ..., F_d$, ucp maps $\psi_i:A\to F_i,$ order zero maps $\phi_i:F_i\to A$ such that the map $\phi:\oplus F_i\to A$ defined by $\phi((y_i)) = \sum \phi_i(y_i)$ is cpc and for $\psi=(\psi_0, \cdots, \psi_d):A\to \oplus F_i$ we have
\begin{enumerate}
\item $\| \psi(f(y)) - f(\psi(y)) \| <\ve$ for $y\in \{ x^*x, xx^* \}.$
\item $\| \phi\circ \psi(y) - y\| <\ve$ for all $y\in \{ f(x^*x), f(xx^*) \}.$
\item $\|  \psi(yy^*)-\psi(y)\psi(y^*) \|<\delta$ for all $y\in \{x, x^*\}.$
\end{enumerate}
In item (1) we are approximating $f$ by polynomials and then applying \cite[Proposition 5.1]{Kirchberg04} to a sufficiently large finite set. For each $i=0, .., d$ obtain a unitary $u_i\in F_i$ such that $u_i\psi_i(x)\psi_i(x)^*u_i^* = \psi_i(x)^*\psi_i(x).$ Such a unitary exists by considering the polar decomposition of $\psi_i(x).$ 

By \cite[Proposition 2.1]{Kirchberg04} for each $i=0,...,d$ there is a *-homomorphism $\pi_i:CF_i\to A$  such that $\phi_i(y) = \pi_i(\id\otimes y)$ for all $y\in F_i.$ For $i=0, ..., d$ define 
\begin{equation*}
z_i=\pi_i(\id^{\frac{1}{2}}\otimes u_if(\psi_i(x)\psi_i(x^*))^{\frac{1}{2}})^*.
\end{equation*}
We use $a \approx_\eta b$ to mean $\| a-b \|<\eta.$ Then
\begin{align*}
z_i^*z_i &= \pi_i(\id\otimes u_if(\psi_i(x)\psi_i(x^*))u_i^*)\\
&=\pi_i(\id\otimes f(\psi_i(x^*)\psi_i(x)))\\
&=\phi_i(f(\psi_i(x^*)\psi_i(x))) \\
&\approx_\ve \phi_i(f(\psi_i(x^*x)))  \\
&\approx_\ve \phi_i(\psi_i(f(x^*x)))
\end{align*}
We then have
\begin{equation*}
\sum_{i=0}^d z_i^*z_i \approx_{2\ve(d+1)} \sum_{i=0}^d \phi_i(\psi_i(f(x^*x))) \approx_\ve f(x^*x).
\end{equation*}
Similarly $\sum_{i=0}^d z_iz_i^* \approx_{2\ve(d+1) + \ve} f(xx^*).$ It is now clear how to choose sequences $z_{i,m}$ satisfying Definition \ref{def:CPdcomp}.
\end{proof}

\begin{definition} \label{def:basicext} Let $A$ be a unital C*-algebra with non-trivial projection $p\in A.$ Let $S\in B(\ell^2(\N))$ be the unilateral shift.  Define $T = S\otimes p + S^*\otimes (1-p),$ and then  $E_{A,p} = C^*(T, \mathcal{K}\otimes A).$  
\end{definition}
\begin{lemma} \label{lem:surjtoep} Let $A, p$ be as in Definition \ref{def:basicext}. Then $E_{A, p}/(\mathcal{K}\otimes A) \cong C(\T).$ 
\end{lemma}
\begin{proof} Let $f\in C(\T)$ and define $f^\vee(z) = f(\overline{z}).$ We let $T_f$ denote the Toeplitz operator with symbol $f.$  By the basic theory of Toeplitz operators, for every $x\in E_{A, p}$ there is a unique $f_x\in C(\T)$ and $y\in \mathcal{K}\otimes A$  such that $x = T_{f_x}\otimes p + T_{f_x^\vee}\otimes (1-p) + y.$  It follows that $x\to f_x$ defines a surjective *-homomorphism.
\end{proof}

\begin{definition} Let $A$ be a unital AF-algebra.  We say $A$ is 2-symmetric if it admits a decomposition with connecting matrices of the form $\left[  \begin{array}{cc} k_n & \ell_n\\ \ell_n & k_n \end{array}  \right]$ for some sequences of integers $k_n, \ell_n.$
\end{definition}

We later provide a K-theoretic and more general proof of Lemma \ref{lem:sqd}.
We prove this case directly because I prefer operator theoretic proofs of quasidiagonality whenever reasonable,  and  the proof highlights how 2-symmetric conspires to produce the quasidiagonal projections.
After the importance of 2-symmetric is realized, which we show in the next lemma,  it is essentially the same proof that $S^{\oplus n}\oplus S^{*\oplus n}$ is quasidiagonal.

\begin{lemma}\label{lem:rankcalc} Let $A = \lim M_{a_n}\oplus M_{a_n}$ be a 2-symmetric AF-algebra with connecting matrices $\left[  \begin{array}{cc} k_n & \ell_n\\ \ell_n & k_n \end{array}  \right].$ Let $\pi_{m,n}:M_{a_m}\oplus M_{a_m}\to M_{a_n}\oplus M_{a_n}$  be the corresponding *-homomorphisms. 
For $x\in M_{a_m}\oplus M_{a_m}$, let $\operatorname{rank}_m(x)$ denote the usual rank.
Let $p \in M_{a_m}\oplus M_{a_m}$ be a projection with $\operatorname{rank}_m(p) = \operatorname{rank}_m(1-p).$  Then for all $n\geq m$ we have $\operatorname{rank}_n(\pi_{m, n}(p)) = \operatorname{rank}_n(\pi_{m, n}(1-p)).$
\end{lemma}
\begin{proof} 
For any $k\geq 0$, let $z_k = (1_{M_{a_k}}, 0)\in  M_{a_k}\oplus M_{a_k}.$ For $x\in M_{a_k}\oplus M_{a_k}$ define  $\operatorname{rank}_{k, 0}(x) = \operatorname{rank}_k(xz_k)$ and  $\operatorname{rank}_{k, 1}(x) = \operatorname{rank}_k(x(1-z_k)).$  

Clearly it is enough to show the conclusion for $n=m+1.$ For any $x\in M_{a_m}\oplus M_{a_m}$ we have
\begin{align*}
\operatorname{rank}_{m+1}(\pi_{m, m+1}(x)) & = \operatorname{rank}_{m+1, 0}(\pi_{m, m+1}(x)) +\operatorname{rank}_{m+1, 1}(\pi_{m, m+1}(x))\\
& = (k_m\operatorname{rank}_{m, 0}(x) +\ell_m\operatorname{rank}_{m, 1}(x))\\
& \hspace{0.2in} + (\ell_m\operatorname{rank}_{m, 0}(x) +k_m\operatorname{rank}_{m, 1}(x))\\
& = (\ell_m+k_m)\operatorname{rank}(x).
\end{align*}
Hence $\operatorname{rank}_{m+1}(\pi_{m, m+1}(p)) = \operatorname{rank}_{m+1}(\pi_{m, m+1}(1-p)).$
\end{proof}

\begin{lemma}\label{lem:sqd} Let $A = \lim M_{a_n}\oplus M_{a_n}$ be a 2-symmetric AF-algebra. Let $k_0 \geq 1$ and set $p = (1_{a_{k_0}}, 0)\in M_{a_{k_0}}\oplus M_{a_{k_0}}\subseteq A.$ 
Then $E_{A, p}$ is quasidiagonal. If $A$ is simple, then $E_{A, p}$ is strongly quasidiagonal. 
\end{lemma}
\begin{proof} Represent $A\subset B(H)$ faithfully with cyclic and separating vector $\eta.$ For each $n> k_0$, let $P_n\in B(H)$ be the finite rank projection onto $\operatorname{span}\{ a\eta: a\in M_{a_n}\oplus M_{a_n}  \}.$  Then the sequence $(P_n)$ approximately commutes with every element of $A.$  Since $\eta$ is cyclic, the sequence $P_n$ converges strongly to the identity. Moreover, since $\eta$ is separating and $A$ is 2-symmetric  we have $r_n:=\operatorname{rank}(pP_n) = \operatorname{rank}((1-p)P_n)$ by Lemma \ref{lem:rankcalc}.

Russell Smucker characterized quasidiagonal weighted shifts \cite{Smucker82} and in the case of quasidiagonality provided explicit projections.  We adapt his construction here.

For each $n,$ let $H_n^+ = \operatorname{span}\{ px\eta:x\in  M_{a_n}\oplus M_{a_n}\}$ and $H_n^- = \operatorname{span}\{ (1-p)x\eta:x\in  M_{a_n}\oplus M_{a_n}\}.$  Fix finite orthonormal bases $\{ \xi_i^+ : i =1,..., r_n \}$
 and $\{\xi_i^- : i =1,..., r_n\}$ for $H_n^+$ and $H_n^-$ respectively. 
 
 For each $1\leq i \leq r_n$ and $k=1, ..., n-1$ define
 \begin{equation*}
 V_{k, i} = \sqrt{\frac{n-k}{n}}\delta_{n+k}\otimes \xi_i^+ + \sqrt{\frac{k}{n}}\delta_{2n-k+1}\otimes \xi_i^-
 \end{equation*}
 
Let $Q_n$ project onto the span of the following orthonormal set 
\begin{equation*}
\mathcal{B} = \{ V_{k, i}: 1\leq i\leq r_n, 1\leq k\leq n-1  \} \cup \{ \delta_k\otimes\xi_i^\pm : 1\leq k\leq n, 1\leq i\leq r_n  \}\cup\{   \delta_{n+1}\otimes\xi_i^-: 1\leq i\leq r_n \}.
\end{equation*}
It is clear that $Q_n$ converges strongly to the identity.

Let $x\in M_{a_n}\oplus M_{a_n}.$ Let $y\in B(\ell^2(\N))$  such that $\ker(y), \ker(y^*)\supseteq \operatorname{span}\{ \delta_k:1\leq k\leq n \}^\perp.$  Then $y\otimes x$ commutes with $Q_n.$
We now show $T = S\otimes p+S^*\otimes (1-p)$ approximately commutes with $Q_n.$ 

We check that each element of $\mathcal{B}$ is approximately returned to $\mathcal{B}.$ We proceed in cases.
\newline
\emph{Case 1:} Let $1\leq i\leq r_n$ and $1\leq k \leq n-2.$ Then
\begin{equation*}
\| TV_{k, i} - V_{k+1, i} \|^2 = \left(  \sqrt{\frac{n-k}{n}} - \sqrt{\frac{n-(k+1)}{n}}  \right)^2 + \left( \sqrt{\frac{k}{n}} - \sqrt{\frac{k+1}{n}}\right)^2 \leq \frac{2}{n}.
\end{equation*}
\emph{Case 2.}  Let $1\leq i\leq r_n$ and $ k = n-1.$ Then
\begin{equation*}
\| TV_{n-1, i} - \delta_{n+1}\otimes \xi_i^- \|^2 = \frac{1}{n} + \left( 1-\sqrt{\frac{n-1}{n}} \right)^2\leq \frac{2}{n}.
\end{equation*}
\emph{Case 3.} Let $1\leq i\leq r_n, 1\leq k\leq n-1, \ve_i\in \{ +, - \},$ but exclude the case $k=1$ and  $\ve_i = -.$
\newline
Then $T(\delta_k\otimes \xi_i^{\ve_i}) = \delta_{k \hspace{0.01in}\ve_i \hspace{0.01in} 1}\otimes \xi_i^{\ve_i}.$
\newline
\emph{Case 4.} Leftover cases. Let $1\leq i\leq r_n$
\begin{enumerate}
\item $T(\delta_k\otimes \xi_i^-) = \delta_{k-1}\otimes \xi_i^-$ for $k = n, n+1.$
\item $T(\delta_1\otimes \xi_i^-)=0.$
\item $\| T(\delta_n\otimes \xi_i^+) - V_{1, i} \|^2 \leq \frac{2}{n}$, similar to Case 2.
\end{enumerate}
By an unilluminating calculation that we omit (see \cite[Lemma B]{Smucker82} for the details) it follows that $\| (1-Q_n)TQ_n \|\leq \sqrt{\frac{2}{n}}.$ Very similar arguments show that $\| (1-Q_n)T^*Q_n \|\leq \sqrt{\frac{2}{n}}.$ Hence $\| [T, Q_n] \|\leq \sqrt{\frac{2}{n}}.$

This shows that $E_{A, p}$ has a faithful quasidiagonal representation.  If $A$ is simple, then any representation of $E_{A, p}$
that doesn't kill $\mathcal{K}\otimes A$ is essential and therefore quasidiagonal by Voiculescu's theorem \cite{Voiculescu76} (see \cite[Theorem 7.2.5]{Brown08}). Any representation that kills $\mathcal{K}\otimes A$ factors through $C(\T)$ and is therefore quasidiagonal.
\end{proof}

Spielberg's \cite[Lemma 1.5]{Spielberg88} shows for an extension $0\to I\to E\to E/I\to 0$ that if both $I$ and $E/I$ are stably finite and $\partial(K_1(E/I))\cap K_0(I)^+ = \{0\}$ then $E$ is stably finite as well.  He also shows that when $E$ satisfies the UCT and either $I$ or $E/I$ is AF then $E$ is in fact quasidiagonal (actually AF-embeddable \cite[Theorem 1.15]{Spielberg88}).

N. Brown and Dadarlat generalized Spielberg's result from AF-ideals  to those that have the $K_0$-embedding property \cite[Theorem 4.6]{Brown04a}.  Recently Moutzouris expanded the known class of C*-algebras with the $K_0$-embedding property.  We recall his definition

\begin{definition}\cite{Moutzouris24}. \label{def:classY} Let $\mathcal{Y}$ be the class of C*-algebras that can be written as a finite direct sum of algebras in the following two classes
\begin{enumerate}
\item $D\rtimes \Z$, where $D$ is separable, nuclear, unital quasidiagonal, satisfies the UCT and the action is minimal.
\item separable ASH algebras.
\end{enumerate}
\end{definition}

\begin{theorem} \label{thm:indexqd} Let $\mathcal{Q}$ be the universal UHF algebra and $A$ be a simple unital infinite dimensional C*-algebra such that $A\otimes \mathcal{Q}$ is locally approximated by algebras in the class $\mathcal{Y}$ (Definition \ref{def:classY}). In particular this covers all AF-algebras $A.$
Let $p\in A$ be a projection such that  $\Z([1_A-p]_0 - [p]_0 )\cap K_0(A)^+ = \{ 0 \}.$ Then $E_{A, p}$ is strongly quasidiagonal.
\end{theorem}
\begin{proof} First note that the index map $\partial:K_1(C(\T))\to K_0(\mathcal{K}\otimes A) $ associated with the extension $0\to \mathcal{K}\otimes A\to E_{A, p}\to C(\T)\to 0$ satisfies
\begin{equation*}
\partial([u]_1) = [e_{11}\otimes (1_A - p)]_0 - [e_{11}\otimes p]_0,
\end{equation*}
where $u = \id_{\T}.$ Hence by assumption and \cite[Lemma 1.5]{Spielberg88}, $E_{A,p}$ is stably finite.  By \cite[Theorem 1.3]{Moutzouris24} $E_{A, p}$ is quasidiagonal.  Since $A$ is unital, infinite-dimensional and simple, the ideal $\mathcal{K}\otimes A$ is essential and hence
strong quasidiagonality follows as in the last paragraph of the proof of Lemma \ref{lem:sqd}.
\end{proof}

\begin{lemma}\label{lem:infdr}  Let $A$ be a unital C*-algebra with a non-trivial projection $p\in A$ and a tracial state $\tau$ such that $\tau(p)\neq 1/2.$  Then $E_{A, p}$ does not satisfy $(d, 1-t)$-comparison for any $d$ and hence has infinite decomposition rank. 
\end{lemma}
\begin{proof}  Let $T\in E_{A, p}$ as in Definition \ref{def:basicext}. Then $1-T^*T= e_{11}\otimes (1-p)$ and $1-TT^*=e_{11}\otimes p$ where $e_{11}$ is the rank one projection onto  $\C\delta_1.$

Suppose that $E_{A, p}$ has $(d, 1-t)$-comparison. By Definition \ref{def:CPdcomp}  with $T=x^*$ there are a sequence of contractions $z_{0, m}, ..., z_{d, m}\in E_{A, p}$ such that
\begin{equation*}
\sum_{i=0}^d z_{i, m}^*z_{i,m} \to e_{11}\otimes p, \qquad \sum_{i=0}^d z_{i, m}z_{i,m}^* \to e_{11}\otimes (1-p).
\end{equation*} 
We can't immediately apply the trace since $z_{i,m}$ are not necessarily in the domain of the densely defined trace on $\mathcal{K}\otimes A$ so we slightly modify them.  Since
\begin{equation*}
\lim_{m\to \infty}(1-e_{11})\otimes 1_A \left(\sum_{i=0}^d z_{i,m}^*z_{i,m}\right)(1-e_{11})\otimes 1_A=0,
\end{equation*}
we have
\begin{equation*}
\lim_{m\to \infty}\left(\sum_{i=0}^d z_{i,m}[(1-e_{11})\otimes 1_A]z_{i,m}^*\right)=0,
\end{equation*}
which implies 
\begin{equation*}
\lim_{m\to \infty}\left(\sum_{i=0}^d (e_{11}\otimes 1_A)z_{i,m}(e_{11}\otimes 1_A)z_{i,m}^*(e_{11}\otimes 1_A)\right)=\lim_{m\to \infty}\sum_{i=0}^d z_{i,m}z_{i,m}^* = e_{11}\otimes (1-p).
\end{equation*}
Similarly,
\begin{equation*}
\lim_{m\to \infty}\left(\sum_{i=0}^d (e_{11}\otimes 1_A)z_{i,m}^*(e_{11}\otimes 1_A)z_{i,m}(e_{11}\otimes 1_A)\right)=\lim_{m\to \infty}\sum_{i=0}^d z_{i,m}^*z_{i,m} = e_{11}\otimes p.
\end{equation*}

Hence setting $y_{i, m} = (e_{11}\otimes 1_A)z_{i, m}(e_{11}\otimes 1_A)$ we obtain 
\begin{equation*}
\sum_{i=0}^d y_{i, m}^*y_{i,m} \to e_{11}\otimes p, \qquad \sum_{i=0}^d y_{i, m}y_{i,m}^* \to e_{11}\otimes (1-p).
\end{equation*} 
Since $y_{i,m}\in e_{11}\otimes A$ we have $\tau(p)=\tau(1-p)$ a contradiction. Hence $E_{A,p}$ does not have $(d, 1-t)$-comparison for any $d$ and hence has infinite decomposition rank by Proposition \ref{prop:basicdr}.
\end{proof}

By Proposition \ref{prop:sr1} every C*-algebra with stable rank 1 has $(0, 1-t)$-comparison.  Hence our method will never produce an example with stable rank 1, so stable rank 2 is the minimal stable rank example we can obtain.

\begin{theorem}\label{thm:main} Let $A = \lim M_{a_n}\oplus M_{a_n}$ be a simple 2-symmetric AF-algebra with more than one tracial state.  Then there is a projection $p\in A$ such that
\begin{enumerate}
\item $\dimnuc E_{A, p} = 1.$
\item $E_{A, p}$ is strongly quasidiagonal.
\item $E_{A, p}$ violates $(d, 1-t)$-comparison  for all $d$, hence $\operatorname{dr} E_{A, p} = \infty.$
\item $E_{A, p}$ has stable rank 2.
\end{enumerate}
\end{theorem}
\begin{proof} We first show how item (1) follows from \cite{Evington26}. Fix any $n\in \N$ and let $p\in M_{a_n}\oplus M_{a_n}$ be a projection (there are no restrictions on $p$ for this part of the argument). Let $m\geq n$ and $I_m =  \mathcal{K}\otimes (M_{a_m} \oplus M_{a_m}).$ Let $E_m = C^*(T, I_m)\subseteq E_{A, p}.$ As in Lemma \ref{lem:surjtoep} for every  $x\in E_m$ there is a unique continuous function $f_x\in C(\T)$ and  $x_I\in I_m$ such that $x = T_{f_x}\otimes p + T_{f_x^\vee}\otimes (1_A-p) + x_I,$ where $f^\vee(z) = f(\overline{z}).$

The homomorphism $x\to f$ provides the short exact sequence
\begin{equation*}
0\to I_m \to E_m \to C(\T)\to 0.
\end{equation*}
Let $e_m = 1_{M_{a_m}}\oplus 0$ and $f\in C(\T)$ non-zero.  Let $\mathcal{C}$ denote the Calkin algebra.
Then the first factor of the Busby invariant $\sigma:C(\T)\to \mathcal{C}\oplus \mathcal{C}$ is given by 
\begin{equation*}
\sigma_1(f) = 1_{\ell^2}\otimes e_m \left(T_f\otimes p + T_{f^\vee}\otimes (1_A-p)\right) + I_m = T_f\otimes e_mp + T_{f^\vee}\otimes e_m(1_A-p) + I_m.
\end{equation*}
Since either $e_mp\neq 0$ or $e_m(1_A-p)\neq 0$ by the simplicity of the Calkin algebra $\sigma_1(f)$ generates the Calkin algebra as an ideal.  Similarly $\sigma_2(f)$ generates the Calkin algebra as an ideal. Therefore the extension is full.
Clearly $I_m$ is a stable C*-algebra hence by \cite[Theorem C]{Evington26} we have $\dimnuc E_m =1$ (since $E_m$ is not AF).  Since $E_{A, p}$ is the inductive limit of $E_m$ it follows that $\dimnuc E_{A, p} = 1$ (since $E_{A, p}$ is not AF).

Since $A$ admits more than one tracial state there is a $k_0$ and a tracial state $\tau$ on $A$ so $\tau((1_{a_{k_0}}, 0))\neq 1/2.$ Set $p = (1_{a_{k_0}}, 0).$
Then Lemma \ref{lem:sqd} shows $E_{A, p}$ is strongly quasidiagonal and 
Lemma \ref{lem:infdr} shows $E_{A, p}$ does not satisfy $(d, 1-t)$-comparison for any $d.$

Finally, by Rieffel's \cite[Theorem 4.11]{Rieffel83a}, the stable rank of $E_{A,p}$ is either 1 or 2, but it can't be 1 by Proposition \ref{prop:sr1}.
\end{proof}

\begin{remark} It seems likely that one could adapt the method of Brake and Winter used for the Toeplitz algebra in \cite{Brake19} to show that $\dimnuc E_{A, p} = 1.$  
\end{remark}

\begin{corollary} \label{cor:main} There are continuum many strongly quasidiagonal C*-algebras $E$  such that
$\dimnuc E = 1$, $E$ violates $(d, 1-t)$-comparison for all $d$ and hence $\operatorname{dr} E=\infty.$
\end{corollary}
\begin{proof} By \cite[Exercise 7.7.4]{Blackadar98} (which comes from Fack and Marechal's  \cite{Fack79}) there are continuum many stable isomorphism classes of simple 2-symmetric AF-algebras with more than one tracial state.
Since the ideal $\mathcal{K}\otimes A \trianglelefteq E_{A,p}$ is the unique simple non-zero ideal it follows that there are continuum many isomorphism classes of $E_{A,p}.$
Hence the conclusion follows from Theorem \ref{thm:main}. 
\end{proof}

\begin{example}\label{ez:easyex} The easiest example that is covered by Corollary \ref{cor:main} is probably the AF-algebra $A$ with dimension group $\Q^2$, order unit $(1, 1)$ and strict ordering in each coordinate.  Then let $p\in A$ be a projection with, for example, $[p]_0 = (1/2, 1/4).$  Then $\Z([1_A-p]_0 - [p]_0 ) = (0, \frac{1}{2}\Z).$ Let $\tau$ be the trace on $A$ such that $\tau_*$ is the second coordinate projection on $\Q^2$ to see $E_{A, p}$ has infinite decomposition rank.
\end{example}

If we are willing to forfeit nuclear dimension 1 we obtain a broader class of examples by Moutzouris's theorem
\begin{corollary}\label{cor:general} Let $A$ be a simple unital C*-algebra with nuclear dimension $d<\infty$ such that $A\otimes \mathcal{Q}$ is locally approximated by algebras in $\mathcal{Y}$ (Definition \ref{def:classY}). Suppose $p\in A$ is a projection such that
\begin{enumerate}
\item $\Z([1_A-p]_0 - [p]_0 )\cap K_0(A)^+ = \{ 0 \}$
\item $\tau(p)\neq \frac{1}{2}$ for some tracial state $\tau$ on $A.$
\end{enumerate}
Then $E_{A,p}$ is strongly quasidiagonal, satisfies $\dimnuc E_{A,p}\leq d+2$, does not have $(r, 1-t)$-comparison for any $r$ and hence has infinite decomposition rank.
\end{corollary}
\begin{proof} This follows from  Theorem \ref{thm:indexqd},  Lemma \ref{lem:infdr} and the nuclear dimension extension bound \cite[Proposition 2.9]{Winter10}
\end{proof}



\end{document}